\def\O{\Omega}

\def \dd#1{{\bf#1}}

\def\cl#1{{\cal#1}}



\def\ouv#1{\smash{\mathop{#1}\limits^{\lower 1pt\hbox
{$\scriptscriptstyle\circ$}}}}

\def\hfl#1#2{\smash{\mathop{\hbox to 12mm{\rightarrowfill}}
\limits^{\scriptstyle#1}_{\scriptstyle#2}}}


\long\def\eno#1#2{\par\smallskip{\bf{#1}}{\it\ {#2}}\par\medskip}

\def\stit#1{\vskip 3mm plus 1mm minus 2mm {\bf{#1}}
		\smallskip}

\font\tir=cmbx10 at 12pt

\def\ref#1#2#3#4{{\bf #1}{\ #2}{\it ,\ #3}{,\ #4}\medskip}


\def \picture #1 by #2 (#3){\midinsert \centerline 
{\vbox to #2{\hrule width #1 heigth 0pt 
depth 0pt \null \vfill \special {picture #3}}}\endinsert }

\def\scaledpicture #1 by #2 (#3 scaled #4) {{
\dimen0 =#1 \dimen1 =$2
\divide \dimen0 by 1000 \multiply \dimen0 by #4
\divide \dimen1 by 1000 \multiply \dimen1 by #4
\picture \dimen0 by \dimen1 (#3 scaled $4)}}

\def\figure #1 #2 #3 {\midinsert \vglue 3mm 
{\vbox to #3 {\hrule width 6cm height 0cm depth 0cm \vfill
{\special {picture #1 scaled #2}}}}\vglue 2mm \endinsert}

\magnification=1200


\overfullrule=0pt

{\centerline {\tir {A note on holomorphic extensions}}}

\bigskip
\bigskip

{\centerline {{\bf R. P\'erez-Marco}\footnote {*} 
{  UCLA, 
Dept. of Mathematics, 
405, Hilgard Ave., Los Angeles, CA-90095-1555, 
USA, e-mail: ricardo@math.ucla.edu; CNRS UMR 8628, Universit\'e 
Paris-Sud, Mathematiques, 91405-Orsay, France.}}

\bigskip
{\bf Abstract.} {\it We give a criterium of holomorphy 
for some type formal power series. This gives
a stronger form of a Rothstein's type extension theorem for  
a particular ring of holomorphic functions.
}

\bigskip


\medskip


\bigskip
\bigskip

We consider the set $R\subset 
\dd C [[z_1, z_2]]$, $z_1, z_2 \in \dd C^k \times 
\dd C^m$ of formal power series of the form
$$
f(z)=f(z_1, z_2)=\sum_n P_n (z_2) z_1^n
$$
where $P_n$ is a polynomial in $m$ variables of total degree
${\hbox {\rm d}}^0 \ P_n \leq C_0+C_1 ||n|| $, for some constants
$C_0 , C_1 >0$.
One easily checks that $R$ is a local sub-ring of $\dd C [[z_1, z_2]]$.
For the notion of $\Gamma$-capacity, that generalizes the 
notion of capacity in one complex variable, we refer to [Ro].

\eno {Theorem.}{Let 
$$
f(z_1, z_2) =\sum_{n} P_n (z_2) z_1^n
$$
be a formal power series of the two complex variables $(z_1, z_2)\in 
\dd C^k\times \dd C^m$. We assume that $(P_n)$ is a sequence of polynomials 
in $m$ variables of total  
degree 
$$
{\hbox {\rm deg}} \ P_n \leq C_0+C_1 ||n|| \ .
$$

We assume that for a set $K\in \dd C^m$ of positive $\Gamma$-capacity, 
$z_2 \in K$ being fixed, the formal power series $f(z_1, z_2)$
converges.

Then for some $C_2 >0$, 
the formal power series $f$ defines a holomorphic function 
in a neighborhood of the axes $\{ z_1=0 \}$ of the form,
$$
U=\{ (z_1, z_2)\in \dd C^k\times \dd C^m ; ||z_1|| \leq {C_2 \over 1+||z_2|| } \} \ .
$$
}

Compare with Rothstein's
theorem (see [Siu] p.25).
Our theorem is motivated and 
has applications in problems of holomorphic dynamics
and small divisors ([PM]) where power series in the ring 
$A$ appear naturally.

\stit {$\Gamma$-capacity.}

We refer to [Ro].
Let $K \subset \dd C^m$. The $\Gamma$-projection of $K$ on 
$\dd C^{m-1}$ is the set $\Gamma_m^{m-1} (K)$ 
of $z=(z_1, \ldots , z_{m-1})\in \dd C^{m-1}$
such that 
$$
K\cap \{ (z, w)\in \dd C^m \}
$$
has positive capacity in the complex plane 
$ \dd C_z=\{ (z, w)\in \dd C^m \}$.
We define 
$$
\Gamma_m^1 (K)=\Gamma_2^1 \circ \Gamma_3^2\circ \ldots \Gamma_m^{m-1} (K) \ .
$$
Finally, the $\Gamma$-capacity is defined as 
$$
\Gamma {\hbox {\rm -Cap}} (K)=\sup_{A\in U(m,\dd C)} {\hbox {\rm Cap}}
\ \Gamma_m^1 (A(K)) \ .
$$
where $A$ runs over all unitary transformations of $\dd C^m$.

We have the following lemma ([Ro] Lemma 2.2.8 p.92)

\eno {Lemma.}{Let $K \subset \dd C^m$, $K\not=\dd C^m$ and assume 
that the intersection of $K$ with any complex line which is not 
a subset of $K$ has inner capacity zero. Then the $\Gamma$-capacity 
of $K$ is zero.}

Thus we are reduced to prove the theorem for $m=1$

\stit {Bernstein lemma.}

We recall (see [Ra] p.156): 

\eno {Lemma (Bernstein).}{Let $K\subset \dd C$ be a non-polar 
set, and $\O$ be the component of ${\overline {\dd C}}-K$ containing
$\infty$.

If $P$ is a polynomial of degree $n$, then for $z\in \dd C$
$$
|P(z)|\leq ||P||_{C^0(K)} e^{ng_{\O} (z, \infty )}
$$
where $g_{\O }$ is the Green function of $\O$.}

\stit {Proof of the theorem.}

We are reduced to prove the theorem for $m=1$.
For $z_2\in K$, let $R(z_2)$ be the radius of convergence
in $z_1$ of $f(z_1,z_2)$.
Let $K_i=\{z_2\in K ; R(z_2) \geq 1/i \}$.
Since a countable union of polar sets is polar, there is 
$K_i$ non-polar. We can take a non-polar sub-compact
$L\subset K_i$ so that there exists $\rho_0 >0$
such that for all $z_2 \in L$
$$
\limsup_{||n||\to +\infty } |P_n (z_2)| \rho_0^{-||n||} <+\infty \ .
$$
Define 
$$
\varphi (z_2) =\limsup_{||n||\to +\infty} |P_n (z_2)| \rho_0^{-||n||} \ .
$$
The function $\varphi$ is lower semi-continuous, and 
$$
L=\bigcup_{p\geq 1} L_p
$$
where $L_p=\{ z\in L ; \varphi (z_2)\leq p \}$ is closed.
By Baire theorem for some $p$, $L_p$ has non-empty 
interior (with respect to $L$), thus some $L_p$ has 
positive capacity. Finally we found a compact set 
$C=L_p$ of positive capacity such that there exists 
$\rho_1>0$ such that for any $z_2\in C$ and $n$,
$$
|P_n(z_2)|\leq \rho_1^{||n||} \ .
$$
Now using Bernstein lemma we conclude that for any 
$z_2 \in \dd C$, for all $n$,
$$
|P_n(z_2)| \leq \rho_1^{||n||} e^{(C_0 +C_1 ||n||) g_{\O} (z_2 ,\infty )} \ .
$$
Finally using the asymptotic
$$
g_{\O } (z_2, \infty )=\log |z_2| +\cl O (1)
$$
we obtain the extension to the desired domain.
\bigskip

\stit {Remark.}

{\bf 1.}
We can improve on the domain of extension 
if we control the growth of the degrees of the 
polynomials $(P_n)$. For instance, the same proof 
shows that if 
$$
\limsup_{n} {1\over ||n||} {\hbox {\rm deg}}\  P_n =0
$$
we have a holomorphic extension to a domain
$$
U=\{ (z_1, z_2)\in \dd C^k \times \dd C^m ; ||z_1||\leq C \}
$$
for some $C>0$.

\medskip

{\bf 2.} As N. Sibony has pointed out recently to me, the 
condition  positive $\Gamma$-capacity in the theorem can be
replaced by non-pluri-polar set  (which is stronger and more 
natural) using the definition of capacity in higher dimension
and the techniques of [Al] and [Si]. Essentially  one 
writes down a general Bernstein lemma in higher dimension
(similar to lemma 6.5 in [Al])
and use it as we do in dimension 1. 
The main dynamical applications, where 
we seek for generic conditions, work as well with the version
with $\Gamma$-capacity. For this first version we content ourselves with 
the above statement.

\bigskip

{\bf Acknowledgments.}

I am grateful to  N. Sibony who told me about 
the notion of $\Gamma$-capacity 
and the reference [Ro], and also explained to me how 
to improve the condition in the theorem as pointed out 
above.

\bigskip

{\centerline {\bf BIBLIOGRAPHY}}

\bigskip
\bigskip

\ref{[Al]}{H. ALEXANDER}{Projective capacity} 
{"Recent developments in several complex variables", J.E. 
Fornaess editor, Princeton Univ. Press, 1981, p. 3-27.}

\ref{[PM]}{R. P\'EREZ MARCO}{Total convergence or general 
divergence in Small Divisors} 
{Preprint, 2000.}

\ref{[Ra]}{T. RANSFORD}{Potential theory in the complex plane}{London
Mathematical Society, Student Texts {\bf 28}, Cambridge University 
Press, 1995.}

\ref{[Ro]}{L.I. RONKIN}{Introduction to the theory 
of entire functions of several variables}{ Translations of
Mathematical monographs, American Mathematical 
Society, {\bf 44}, 1974.}

\ref{[Si]}{N. SIBONY}{Sur la fronti\`ere de Sylow des domaines de 
$\dd C^n$}
{Math. Ann., {\bf 273}, 1985, p. 115-121.}

\ref{[Siu]}{Y.-T. SIU}{Techniques of extension of analytic objects}
{Lecture notes in pure and applied mathematics, {\bf 8}, 
Dekker, 1974.}

\end